\documentclass[10pt, letterpaper, conference]{ieeeconf}

\usepackage{amsmath}
\usepackage{amssymb, latexsym, mathrsfs}
\usepackage{fancybox, graphicx, subfigure}
\usepackage{enumerate, endnotes}
\usepackage[usenames, dvipsnames]{color}

\definecolor{myred}{rgb}{0.6, 0, 0}
\definecolor{mygreen}{rgb}{0, 0.5, 0}
\definecolor{myblue}{rgb}{0, 0, 0.5}
\definecolor{mycyan}{rgb}{0, 0.5, 0.5}

\newcommand{\R}{\ensuremath{\mathbb{R}}}

\newcommand{\N}{\ensuremath{\mathbb{N}}}
\newcommand{\posR}{\ensuremath{\R_{\ge 0}}}

\newcommand{\PSet}{\ensuremath{\mathcal{P}}}
\newcommand{\CSet}{\ensuremath{\mathcal{C}}}

\newcommand{\ra}{\ensuremath{\rightarrow}}

\newcommand{\lra}{\ensuremath{\longrightarrow}}
\newcommand{\eps}{\ensuremath{\varepsilon}}
\newcommand{\fa}{\ensuremath{\forall\,}}
\renewcommand{\le}{\ensuremath{\leqslant}}
\renewcommand{\ge}{\ensuremath{\geqslant}}

\newcommand{\therex}{\ensuremath{\exists\,}}
\newcommand{\setmin}{\ensuremath{\!\smallsetminus\!}}

\newcommand{\ClassK}{\ensuremath{\mathcal{K}}}
\newcommand{\ClassKinfty}{\ensuremath{\mathcal{K}_{\infty}}}
\newcommand{\ClassKL}{\ensuremath{\mathcal{KL}}}

\newcommand{\epower}[1]{\ensuremath{\,\mathrm{e}^{#1}}}
\newcommand{\norm}[1]{\ensuremath{\left\lVert #1 \right\rVert}}
\newcommand{\Expec}[1]{\ensuremath{\mathsf{E}\!\left[\vphantom{\big|}#1\vphantom{\big|}\right]}}
\newcommand{\CExpec}[2]{\ensuremath{\mathsf{E}^{#2}_{\vphantom{T}}\!\left[\vphantom{\big|}#1\vphantom{\big|}\right]}}
\newcommand{\CProb}[2]{\ensuremath{\mathsf{P}^{#2}_{\vphantom{T}}\!\left(\vphantom{\big|}#1\vphantom{\big|}\right)}}
\newcommand{\Prob}[1]{\ensuremath{\mathsf{P}\!\left(\vphantom{\big|}#1\vphantom{\big|}\right)}}
\newcommand{\LieD}[2]{\ensuremath{\mathrm L_{#1}{#2}}}
\newcommand{\indic}[1]{\ensuremath{\boldsymbol{1}_{#1}}}

\newcommand{\mc}[1]{\ensuremath{\mathcal{#1}}}

\newcommand{\Lp}[1]{\ensuremath{\boldsymbol{L}_{#1}}}
\newcommand{\secref}[1]{\S\ref{#1}}

\newcommand{\transp}{\ensuremath{^{\scriptscriptstyle{\mathrm T}}}}
\newcommand{\sigalg}{\ensuremath{\mathfrak{F}}}
\newcommand{\cardP}{\ensuremath{\mathrm{N}}}

\newcommand{\xz}{\ensuremath{x_0}}

\newcommand{\drv}{\ensuremath{\,\mathrm{d}}}
\newcommand{\ol}{\overline}

\newcommand{\wt}{\widetilde}

\newcommand{\cadlag}{c{\`a}dl{\`a}g}
\renewcommand{\subset}{\ensuremath{\subseteq}}

\newcommand{\mx}{\ensuremath{\vee}}
\newcommand{\mn}{\ensuremath{\wedge}}

\newcommand{\DefEnd}{\hspace{\stretch{1}}{$\Diamond$}}
\newcommand{\AssumptionEnd}{\hspace{\stretch{1}}{$\diamondsuit$}}

\newtheorem{theorem}{Theorem}

\newtheorem{proposition}[theorem]{Proposition}
\newtheorem{lemma}[theorem]{Lemma}
\newtheorem{defn}[theorem]{Definition}

\newtheorem{assumption}[theorem]{Assumption}

\IEEEoverridecommandlockouts

\newcommand{\iss}{{\sc iss}}

\title{\Large \bf Towards ISS Disturbance Attenuation for Randomly Switched Systems}
\author{
	Debasish Chatterjee and Daniel Liberzon\thanks{The authors are with the Coordinated Science Laboratory, University of Illinois at Urbana-Champaign, USA}
	\thanks{emails: {\tt dchatter@uiuc.edu, liberzon@uiuc.edu}}
	\thanks{This work was supported by NSF's CSR program (Embedded and Hybrid Systems area) under grant NSF-CNS-0614993.}
}

\begin{document}
	\maketitle

	\begin{abstract}
		We are concerned with input-to-state stability (\iss{}) of randomly switched systems. We provide preliminary results dealing with sufficient conditions for stochastic versions of \iss{} for randomly switched systems without control inputs, and with the aid of universal formulae we design controllers for \iss{}-disturbance attenuation when control inputs are present. Two types of switching signals are considered: the first is characterized by a statistically slow-switching condition, and the second by a class of semi-Markov processes.
	\end{abstract}


	\section{Introduction}
		\PARstart{S}{ince} its introduction in~\cite{ref:sontagISS} the concept of input-to-state stability (\iss{}) has received widespread attention on both theoretical and practical fronts; see~\cite{ref:sontagISSnote} for a recent detailed discussion. The \iss{} property characterizes behavior of the state trajectory of a deterministic nonlinear system perturbed by bounded disturbance inputs; as such it provides a framework for robustness analysis of nonlinear systems. Initially stated for deterministic inputs, various extensions of the \iss{} property have been made for inputs modeled as random processes, one of which is exponential input-to-state stability~\cite{ref:tsiniasExpISS}. The \iss{} property has been employed in constructive ways for stability analysis, stabilizing feedback controller synthesis, adaptive control schemes, etc.

		With the growing interest in the theory and applications of hybrid systems, considerable effort has been directed towards understanding the behavior of switched systems. A switched system has two ingredients: a family of subsystems, and a switching signal which specifies the active subsystem at each instant of time. An important control-theoretic issue is that of stability and stabilization of these systems, and a number of interesting techniques have evolved over the past two decades to deal with this; for a discussion see, e.g.,~\cite[Chapters~2, 3]{ref:liberzonbk}.
		More recently, looking beyond stability, robustness and \iss{} properties of deterministic switched systems have received attention; see~\cite{ref:linhAutomatica06} and the references therein. There appears to be a common thread of slow switching in these results. That is to say, if the constituent subsystems are each \iss{} and the switching is sufficiently slow, then the switched system is also \iss{}.

		In this article we are concerned with \iss{} of randomly switched nonlinear systems, i.e., \iss{} properties of switched systems whose switching signal is a random process. We provide preliminary results dealing with sufficient conditions for a stochastic version of \iss{} of these systems. Two types of switching signals are considered; the first is characterized by a statistically slow switching condition, and the second is a class of semi-Markov processes. For these classes of switching signals it is difficult to apply traditional approaches which rely on an infinitesimal (or extended) generator~\cite{ref:davisMarkovModelsOptbk}, for either there is little information available about the parameters of the switching signal, or there is strong dependence on past history of the process.
		
		The approach pursued here employs multiple \iss{}-Lyapunov functions in the spirit of our earlier works~\cite{ref:ransw, ref:stabranswcorr} on stability analysis of randomly switched systems without inputs. Our approach highlights the interaction of deterministic dynamical systems with the stochastic switching signal. The switching signals considered here are adopted from these articles, but the analysis in the presence of inputs as we carry out here is more involved.

		With the analysis results in hand, we turn to control synthesis. Two types of controller architectures are considered: in the first case the controllers depend on both the switching signal and the state, and in the second case the controllers depend only on the state. The technical tools are off-the-shelf universal formulae for \iss{} disturbance attenuation~\cite{ref:liberzonISSclf} and our analysis results. To the best of our knowledge this is the first time that \iss{} under random switching is being studied.

		The article is organized as follows. In~\secref{s:prelims} we fix notations and define our property of interest. The analysis results are given in~\secref{s:ana}, a proof of one result is sketched in~\secref{s:proofs}, and the synthesis results are presented in~\secref{s:syn}. We conclude in~\secref{s:concl} with a short discussion of the case of Markovian switching signals.

	\section{Preliminaries}
		\label{s:prelims}
		Let $\posR := [0, \infty[$, $\norm{\cdot}$ denote the Euclidean norm on $\R^n$, and $\norm{f}_{A}$ denote the essential supremum norm of the function $f$ on the set $A\subset\posR$. Recall that a function $\alpha:\posR\lra\posR$ belongs to class-$\ClassK$ if $\alpha$ is monotone increasing, continuous, and $\alpha(0) = 0$. Also, $\alpha$ belongs to class-$\ClassKinfty$ if $\alpha\in\ClassK$, and $\alpha\nearrow\infty$. A function $\beta:\posR^2\lra\posR$ belongs to class-$\ClassKL$ if $\beta(\cdot, s)\in\ClassK$ for each $s$ and $\beta(r, \cdot)\searrow 0$ for each $r$. We let $x\mn y := \min\{x, y\}$ and $x\mx y := \max\{x, y\}$ for $x, y\in\R$.

		Let $(\Omega, \sigalg, \mathsf P)$ be a probability space~\cite{ref:raoStochProcGenTheo}, with $\Omega$ the set of events, $\sigalg$ a sigma-algebra on $\Omega$, and $\mathsf P$ a probability measure on $(\Omega, \sigalg)$. We let $\Expec{\cdot}$ denote mathematical expectation and $\CExpec{\vphantom{\big|}\cdot\vphantom{\big|}}{\sigalg'}$ (or $\Expec{\cdot\,\big|\,\sigalg'}$) denote conditional mathematical expectation given a sigma-subalgebra $\sigalg'$ of $\sigalg$. We let $\CProb{\cdot}{\sigalg'}$ (or $\Prob{\cdot\big|\sigalg'}$) denote conditional probability given $\sigalg'$.

		\subsection{Randomly switched systems with disturbance inputs}
		Let $\PSet := \{1, \ldots, \cardP\}$ be a finite index set, and for each $i\in\PSet$ let us consider the system
		\begin{equation}
			\dot x = f_i(x, d)
			\label{e:ithsys}
		\end{equation}
		where $f_i:\R^n\times\R^k\lra\R^n$ is a continuously differentiable vector field, $f_i(0, 0) = 0$. We allow $d:\posR\lra\R^k$ to be a measurable, locally essentially bounded function of time; this ensures local existence and uniqueness of solutions of~\eqref{e:ithsys}. Let $\sigma$ be a \cadlag{} stochastic process (i.e., a stochastic process whose sample paths are continuous from the right and possess limits from the left) on $(\Omega, \sigalg)$ taking values in $\PSet$. We assume that for each $t \ge 0$ and each $\omega\in\Omega$ there exists a strictly positive number $\epsilon(t, \omega)$ such that $\sigma(t+s, \omega) = \sigma(t, \omega)$ on $[t, t+\epsilon(t, \omega)[$. Under this condition we know~\cite[Theorem~T26, p.~304]{ref:bremaudPointProc} that the filtration $(\sigalg_t)_{t\ge 0}$ generated by $\sigma$ is right-continuous, and we augment $\sigalg_0$ with all $\mathsf P$-null sets. We say that $\sigma$ is a random switching signal, and it generates the randomly switched system from the family~\eqref{e:ithsys} given by
		\begin{equation}
			\dot x = f_\sigma(x, d), \qquad (x(0), \sigma(0)) = (\xz, \sigma_0),\quad t\ge 0,
			\label{e:ssysdef}
		\end{equation}
		where $\xz\in\R^n$.

		\subsection{Input-to-state stability}
		Input-to-state stability (\iss{}) was formulated for a single system in~\cite{ref:sontagISS}; let us state the definition corresponding to the $i$-th member of the family defined above.

		The system~\eqref{e:ithsys} is \emph{input-to-state stable} if there exist functions $\beta_i\in\ClassKL$ and $\gamma_i\in\ClassKinfty$ such that for every $\xz\in\R^n$ and measurable and locally bounded input $d$, the estimate
		\begin{equation}
			\norm{x(t)} \le \beta_i(\norm{\xz}, t) + \gamma_i\bigl(\norm{d}_{[0, t[}\bigr)
			\label{e:issdef}
		\end{equation}
		holds for all $t\ge 0$ along solutions of~\eqref{e:ithsys}.

		\begin{defn}
			\label{d:issm}
			The system~\eqref{e:ssysdef} \emph{satisfies an \iss{} in $\Lp 1$ estimate at switching instants} if there exist functions $\beta\in\ClassKL$ and $\alpha, \gamma\in\ClassKinfty$ such that for every $\xz\in\R^n$, every measurable and essentially bounded input $d$, the estimate
			\begin{equation}
				\Expec{\alpha(\norm{x(\tau_\nu)})} \le \beta(\norm{\xz}, \nu) + \gamma\bigl(\norm{d}_{\posR}\bigr)
				\label{e:rswiss}
			\end{equation}
			holds for all $\nu\in\N$ along solutions of~\eqref{e:ssysdef}.\DefEnd{}
		\end{defn}

		Notice that the expectation on the left hand side involves a class-$\ClassKinfty$ function $\alpha$. In the absence of randomness this statement in terms of $\alpha$ is equivalent to the statement that $\norm{x(t)} \le \beta'(\norm{\xz}, t) + \gamma'\bigl(\norm{d}_{[0, t[}\bigr)$ for some functions $\beta'\in\ClassKL$ and $\gamma'\in\ClassKinfty$, where we have employed a weak triangle inequality for class-$\ClassKinfty$ functions.\footnote{The weak triangle inequality for a function $\gamma\in\ClassKinfty$ is: $\gamma(r_1 + r_2) \le \gamma(2(r_1\mx r_2)) \le \gamma(2r_1) + \gamma(2r_2)$.} In the context of randomly switched systems, however, without further assumptions on $\alpha$ one cannot conclude that $\Expec{\norm{x(t)}} \le \beta'(\norm{\xz}, t) + \gamma'\bigl(\norm{d}_{[0, t[}\bigr)$ from~\eqref{e:rswiss}. However, it is often the case that we get polynomial functions of the state inside the expectation, which in general yields stronger bounds. For instance, if the function $\alpha$ is quadratic, it is convex, and an application of Jensen's inequality\footnote{Jensen's inequality~\cite{ref:borovkovbk} states that if $X$ is an integrable random variable and $\phi:\R\lra\R$ is a convex function, then $\phi\bigl(\Expec{X}\bigr) \le \Expec{\phi(X)}$.} leads to the last inequality.

		Let us also note that Definition~\ref{d:issm} does not claim \iss{} of every sample path of the system~\eqref{e:ssysdef}; the qualitative and quantitative aspects of this definition do not concern information about individual trajectories.

		Suppose that~\eqref{e:ithsys} is \iss{} for each $i\in\PSet$. Then by definition there exist functions $\beta_i\in\ClassKL$ and $\gamma_i\in\ClassKinfty$ such that~\eqref{e:issdef} holds along solutions of the $i$-th subsystem. However, without further stipulations on $\sigma$, in general it is not true that the switched system generated by $\sigma$ from the family $\{f_i\}_{i\in\PSet}$ retains the \iss{} property (i.e., there will exist \emph{unique} functions $\beta\in\ClassKL$ and $\gamma\in\ClassKinfty$ such that~\eqref{e:rswiss} holds for any trajectory of the switched system~\eqref{e:ssysdef}). In~\secref{s:ana} we consider different classes of switching signals for which we give sufficient conditions for different types of \iss{}-type estimates.

	\section{Analysis Results}
		\label{s:ana}

		\subsection{Statistically slow switching}
		We assume no more structure of the switching signal than a slow switching condition, which is reminiscent of the switching rate of a Poisson counter. A similar condition was employed in the main theorem of~\cite{ref:ransw}, where we dealt with stability under no disturbance inputs and slow switching. We also assume that each member of the family of subsystems is \iss{}. First a piece of notation: let $N_\sigma(t_2, t_1)$ denote the number of jumps made by $\sigma$ on the interval $]t_1, t_2]\subset\posR$, $t_1 \le t_2$.

		Recall that we have $\bigl(\Omega, \sigalg, (\sigalg_t)_{t\ge 0}, \mathsf P\bigr)$ as a complete filtered probability space satisfying the usual conditions. A $[0, \infty]$-valued random variable $\tau$ is an $(\sigalg_t)_{t\ge 0}$-\emph{stopping time} if $\{\tau\le t\}\in \sigalg_t$ for each $t\ge 0$. A random variable $\tau'$ is an $(\sigalg_t)_{t\ge 0}$-\emph{optional time} if $\{\tau' < t\}\in\sigalg_t$ for each $t > 0$. It is quite clear that an $(\sigalg_t)_{t\ge 0}$-optional time is a $(\sigalg_t)_{t\ge 0}$-stopping time. It is a standard result that if $(\sigalg_t)_{t\ge 0}$ is a right-continuous filtration, every $(\sigalg_t)_{t\ge 0}$-stopping time is also an $(\sigalg_t)_{t\ge 0}$-optional time. For details see, e.g.,~\cite{ref:raoStochProcGenTheo}.

		\begin{defn}
			\label{d:sigmagen}
			The switching signal $\sigma$ is said to belong to \emph{class G} if the following condition holds: there exist $\ol\lambda, \wt\lambda > 0$ and $ k_0\in\N\cup\{0\}$, such that for every $(\sigalg_t)_{t\ge 0}$-stopping time $t'$ and for all $k\ge 0$:\\
					$\displaystyle{\CProb{N_\sigma(t'+s, t') = k}{\sigalg_{t'}} \le \epower{-\wt\lambda s}\frac{\bigl(\ol\lambda s\bigr)^k}{k!}}$.\DefEnd
		\end{defn}

		Note that if $\ol\lambda = \wt\lambda$ and $k_0 = 0$, then Definition~\ref{d:sigmagen} gives the jump rate of a stationary Poisson process.

		We have the following
		\begin{lemma}
			If $\sigma$ belongs to class G, then $(\tau_i)_{i\in\N}$ is almost surely divergent.
			\label{l:nozeno}
		\end{lemma}

		One can prove this by estimating the expected value of $N_\sigma(t, 0)$ from the bound in Definition~\ref{d:sigmagen} for a fixed $t \ge 0$ (since each fixed $t$ is an $(\sigalg_t)_{t\ge 0}$-optional time), which is readily seen to be finite. See also~\cite[Chapter~3]{ref:myphdthesis} for an alternative argument.

		Our results employ a family of \iss{}-Lyapunov functions; the following assumption collects the properties we require from them. The analysis will proceed with the aid of \iss{}-Lyapunov-like functions.\footnote{Notice that since we do not always require, (as in~\secref{s:semiMarkov}) every subsystem to be \iss{}, these functions are not \iss{}-Lyapunov functions in the strict sense of the term.} The following assumption collects the properties we shall require from them.
		\begin{assumption}
			Suppose that there exist continuously differentiable functions $V_i:\R^n\lra\posR$, $i\in\PSet$, functions $\alpha_1, \alpha_2, \chi\in\ClassKinfty$, and numbers $\mu\ge 1$, $\lambda_i\in\Lambda\subset\R$, $i\in\PSet$, such that for all $(i, x, d)\in\PSet\times\R^n\times\R^k$ we have
			\begin{enumerate}[{\rm (Vd1)}]
				\item $\alpha_1(\norm x) \le V_i(x) \le \alpha_2(\norm x)$;
				\item $\displaystyle{\frac{\partial V_i}{\partial x}(x) f_i(x, d) \le -\lambda_i V_i(x) + \chi(\norm d)}$;
				\item $V_i(x) \le \mu V_j(x)$.\AssumptionEnd
			\end{enumerate}
			\label{a:Vi}
		\end{assumption}

		Note that if we allow $\Lambda$ to include negative numbers, then not all $\lambda_i$'s need to be positive, which in turn means that not all subsystems are required to be \iss{}.

		The function $V_i$ in (Vd1) and (Vd2) above is called an \emph{\iss{}-Lyapunov function} for the $i$-th subsystem. If $\Lambda$ consists of positive real numbers, (Vd2) is equivalent to each subsystem being \iss{}. Let us note that conventionally \iss{}-Lyapunov functions are defined in a little different way, for instance, the right-hand side of (VL2) is $-\alpha'(\norm x) + \chi'(\norm d)$, or the right-hand side of (VL2) is $-\alpha'(\norm x)$, for $\alpha', \chi'\in\ClassKinfty$, but they are equivalent to (VL2), as proved in~\cite{ref:pralywangISS}.

		\begin{theorem}
			\label{t:issm}
			Consider the switched system~\eqref{e:ssysdef}, and suppose that
			\begin{enumerate}[{\rm (G1)}]
				\item $\sigma$ belongs to class G;
				\item Assumption~\ref{a:Vi} holds with $\Lambda = \{\lambda_\circ\}$, $\lambda_\circ > 0$;
				\item $\mu < \bigl(\wt\lambda + \lambda_\circ\bigr)/\ol\lambda$.
			\end{enumerate}
			Then there exists a monotonically nondecreasing sequence $(T_i)_{i\in\N}$ of $(\sigalg_t)_{t\ge 0}$-optional times with $\lim_{t\ra\infty} T_i = \infty$ a.s., and functions $\beta\in\ClassKL$, $\alpha,\gamma\in\ClassKinfty$, such that
			\begin{align}
				\label{e:star}
				& \Expec{\alpha(\norm{x(t)})\indic{\{t\in[T_{i-1}, T_i[\}\cap\{T_{i-1} < \infty\}}}\nonumber\\
				& \qquad\qquad\qquad\qquad \le \beta(\norm{\xz}, t) \mx \gamma\bigl(\norm{d}_{\posR}\bigr)
			\end{align}
			for all $t\ge 0$ and $i\in\N$.
		\end{theorem}

		The proof of Theorem~\ref{t:issm} is rather long, and may be found in~\cite[Chapter~3]{ref:myphdthesis}; we sketch the main steps in~\S\ref{s:proofs}. See also \S\ref{s:discussion} below for a discussion.
		\subsection{A class of semi-Markov switching signals}
		\label{s:semiMarkov}
		In this subsection we assume $\sigma$ possesses more structure than being statistically slow-switching. Let $S_i := \tau_i - \tau_{i-1}$ for $i\in\N$ be the $i$-th holding time, $(\tau_i)_{i\in\N}$ being the sequence of switching instants.
		\begin{defn}\mbox{}
			\label{d:sigmai}
			The switching signal $\sigma$ is said to belong to \emph{class UH} if it satisfies:
			\begin{enumerate}[{\rm (UH1)}]
				\item the sequence $(S_i)_{i\in\N}$ of holding times is a collection of i.i.d uniform-$(T)$ random variables;
				\item the sequence $(\sigma(\tau_i))_{i\in\N\cup\{0\}}$ of values is i.i.d with $\Prob{\sigma(\tau_i) = i} = q_j$ for some $q_j\in\,]0, 1[$, $j\in\PSet$;
				\item the two sequences $(S_i)_{i\in\N}$ and $(\sigma(\tau_i))_{i\in\N}$ are mutually independent.\DefEnd
			\end{enumerate}
		\end{defn}

		The class UH of switching signals is simply a representative example of the class of semi-Markov switching signals that we can treat in our framework; see~\cite{ref:myphdthesis} for other classes of switching signals and related discussion.

		\begin{lemma}
			For switching signals of class UH, the sequence $(\tau_i)_{i\in\N}$ is almost surely divergent.
			\label{l:nozenoi}
		\end{lemma}

		The above lemma can be established by appealing to the Strong Law of Large Numbers~\cite[Chapter~2]{ref:raoProbTheo}; see also~\cite[Chapter~2]{ref:myphdthesis} for alternative arguments.
		\begin{theorem}
			Consider the switched system~\eqref{e:ssysdef}. Suppose that
			\begin{enumerate}[{\rm (U1)}]
				\item $\sigma$ belongs to class UH;
				\item Assumption~\ref{a:Vi} holds with $\Lambda = \R$;
				\item $\displaystyle{\sum_{j\in\PSet}\frac{\mu q_j\bigl(1-\epower{-\lambda_j T}\bigr)}{\lambda_j T} < 1}$.
			\end{enumerate}
			Then~\eqref{e:ssysdef} satisfies an \iss{} in $\Lp 1$ estimate at switching instants.
			\label{t:issmunif}
		\end{theorem}

	\subsection{Discussion}
		\label{s:discussion}
		The results above fall short of being satisfactory. Indeed, perhaps the most natural adaptation of the \iss{} concept to the stochastic case would involve bounds of the type
		\begin{equation}
			\label{e:ideal}
			\Expec{\alpha(\norm{x(t)})} \le \beta(\norm{\xz}, t) + \gamma\bigl(\norm{d}_{\posR}\bigr)
		\end{equation}
		for all $\xz\in\R^n$, $t\ge 0$, and essentially bounded inputs $d$. However, the technical difficulties, particularly in the absence of Markovian assumptions on $\sigma$, are formidable. Let us consider switching signals belonging to class G. If $B$ denotes the ball around the origin whose radius is $\rho\bigl(\norm{d}_{\posR}\bigr)$, and $B'$ is a larger concentric ball, then the solution trajectory $x(\cdot)$ enters $B$ and exits $B'$ at random instants, as defined in~\eqref{e:tidef}; the sequence $(T_i)_{i\in\N}$ in~\eqref{e:star} is actually this set of random instants. There is no further structure which prevents the number of exit/entry times from increasing at least linearly with time $t$ (the linearity follows at once from the observation that the set of vector fields $\{f_i\}_{i\in\PSet}$ is locally Lipschitz, and that $\norm{d}_{\posR} < \infty$). It is also clear that estimates for the probability distribution of the holding times are not available. Hence ``gain-margin'' type arguments appear to be the only mode of attack, as we pursue in~\secref{s:proofs}. As asserted in Theorem~\ref{t:issm}, it is possible to get bounds on the expectation of the state at some given time, restricted to each of these random excursion intervals, but gluing these estimates to get a uniform bound for a given time $t$ is a difficult problem, and in our case it is yet unsolved.

		On the other hand, in the case of switching signals of class UH, the holding times are explicitly characterized, but the chief issue is that of obtaining an estimate for $\Expec{\alpha(\norm{x(t)})}$ from an \iss{} estimate in $\Lp 1$ at switching instants. To wit, there can potentially be indefinitely many jumps of $\sigma$ before and after a given time $t$; therefore countably many simultaneous interpolations are needed to get an estimate of $\Expec{\alpha(\norm{x(t)})}$, and such an interpolation is again a difficult problem. Unlike in the deterministic case, one is necessarily forced to work with random intervals.

		Let us also note that \iss{}-type estimates ``in probability'' for diffusion processes have appeared in the literature, for instance, in~\cite[Theorem~4.2]{ref:krsticdengbk}, and more recently in~\cite[\S2]{ref:smallgainwrong}. Although the system models in the above references differ from ours, the essential technical difficulties remain the same. Unfortunately, these difficulties were not realized in the aforesaid references, and the claims made in both of them are still open.

	\section{Proofs}
	\label{s:proofs}
		\emph{Proof of Theorem~\ref{t:issm} (Sketch).}
		The argument is divided into five steps for convenience. We shall employ the equivalent ``gain-margin'' characterization~\cite{ref:sontagISSnote} of \iss{} of the individual subsystems; see~\cite[chapter~3]{ref:myphdthesis} for a more detailed proof.

			\textsl{Step~1.} Let us fix an essentially bounded disturbance input signal $d$ with $\norm{d}_{\posR} > 0$, an initial condition $\xz\in\R^n$, and define the open sets $C_1 := \bigl\{z\in\R^n\,\big|\,\norm z < \rho\bigl(\norm{d}_{\posR}\bigr)\bigr\}$ and $C_2 := \bigl\{z\in\R^n\,\big|\,\norm{z} < \eta\rho\bigl(\norm{d}_{\posR}\bigr)\bigr\}$, where $\eta > 0$ is chosen such that $\alpha_1\bigl(\eta\rho\bigl(\norm{d}_{\posR}\bigr)\bigr) > 2\alpha_2\bigl(\rho\bigl(\norm{d}_{\posR}\bigr)\bigr)$. Let us suppose that $\xz\not\in C_1$, the other case being similar. We define the following sequence of random times taking values in $[0, \infty]$:
			\begin{equation}
				\begin{aligned}
				\check t_1 & := \inf\{t > 0\mid x(t)\in C_1\},\\
				\hat t_1 & := \inf\{t > \check t_1\mid x(t)\in \R^n\setmin C_2\},\\
				\ldots\\
				\check t_{i+1} & := \inf\{t > \hat t_{i}\mid x(t)\in C_1\}\quad \text{for }i\in\N,\\
				\hat t_{i+1} & := \inf\{t > \check t_{i+1}\mid x(t)\in \R^n\setmin C_2\}\quad \text{for }i\in\N,
				\end{aligned}
				\label{e:tidef}
			\end{equation}
			where it is understood that if any $\check t_i$ or $\hat t_i$ is $\infty$, then each of the definitions which follow it in the above sequence is set to $\infty$. We note that both $\check t_i$ and $\hat t_i$ are $[0, \infty]$-valued $(\sigalg_t)_{t\ge 0}$-optional times.

			\textsl{Step~2.} Pointwise on $\bigl\{t, \tau_i\in [0, \check t_1[\bigr\}$ we have $x(t), x(\tau_i)\in\R^n\setmin C_1$, and from (Vd2)-(Vd3) we get
			\begin{align*}
				\Expec{V_{\sigma(t)}(x(t))\indic{\{t\in[0, \check t_1[\}}} = \alpha_2(\norm{\xz})\epower{-\bigl(\lambda_\circ + \wt\lambda - \mu\ol\lambda\bigr)t}.
			\end{align*}
			Therefore,
			\[
			\Expec{V_{\sigma(t)}(x(t))\indic{\{t\in[0, \check t_1[\}}} \le \beta(\norm{\xz}, t)\quad \fa t\ge 0,
			\]
			where $\beta(r, s) := \alpha_2(r)\epower{-\lambda s}$, $\lambda := \lambda_\circ + \wt\lambda - \mu\ol\lambda > 0$ by (Gd3).

			\textsl{Step~3.} Pointwise on $\bigl\{t, \tau_i\in[\check t_{j}, \hat t_{j}[\bigr\}\cap\bigl\{\check t_j < \infty\bigr\}$ for $i, j\in\N$ we have $x(t), x(\tau_i)\in C_2$ by~\eqref{e:tidef} and continuity of $x(\cdot)$. Employing (Vd1) leads to
			\[
			\fa t\in[\check t_j, \hat t_j[\quad V_{\sigma(t)}(x(t)) \le \alpha_2\bigl(\eta\rho\bigl(\norm{d}_{\posR}\bigr)\bigr).
			\]
			whenever $\hat t_j < \infty$. Taking expectations we arrive at
			\begin{align*}
				& \Expec{V_{\sigma(t)}(x(t))\indic{\{\hat t_j < \infty\}\cap\{t\in[\check t_j, \hat t_j[\}}}\\
				& \le \alpha_2\bigl(\eta\rho\bigl(\norm{d}_{\posR}\bigr)\bigr)\Prob{\{t\in[\check t_j, \hat t_j[\}\cap\{\hat t_j < \infty\}}.
			\end{align*}

			\textsl{Step~4.} Pointwise on $\bigl\{t, \tau_i\in[\hat t_j, \check t_{j+1}[\bigr\}\cap\bigl\{\hat t_j < \infty\bigr\}$ for $i, j\in\N$ we have 
			\begin{equation}
				\begin{aligned}
				\frac{\partial V_{\sigma(t)}}{\partial x}(x(t))f_{\sigma(t)}(x(t), d(t)) & \le -\lambda_\circ V_{\sigma(t)}(x(t)),\\
				\fa k\in\PSet\quad V_{\sigma(\tau_i)}(x(\tau_i)) & \le \mu V_k(x(\tau_i))
				\end{aligned}
				\label{e:jthexcur}
			\end{equation}
			in view of (Vd2)-(Vd3). Therefore,
			\begin{align}
				& \Expec{V_{\sigma(t)}(x(t))\indic{\{t\in[\hat t_j, \check t_{j+1}[\}\cap\{\hat t_j < \infty\}}}\nonumber\\
				& \le \Expec{\sup_{s\ge 0}V_{\sigma(\hat t_j + s)}(x(\hat t_j + s))\indic{\{\hat t_j + s < \check t_{j+1}\}\cap\{\hat t_j < \infty\}}}.
				\label{e:step40}
			\end{align}
			It can be shown that the process $\bigl(V_{\sigma(\hat t_j + s)}(x(\hat t_j + s))\indic{\{\hat t_j + s < \check t_{j+1}\}\cap\{\hat t_j < \infty\}}\bigr)_{s\ge 0}$ is a nonnegative $\bigl(\sigalg_{\hat t_j + s}\bigr)_{s\ge 0}$-potential. Further detailed calculations lead to
			\[
			\Expec{V_{\sigma(t)}(x(t))\indic{\{t\in[\hat t_j, \check t_{j+1}[\}\cap\{\hat t_j < \infty\}}} \le \gamma\bigl(\norm{d}_{\posR}\bigr),
			\]
			where we let $\gamma(r) := (1+1/\delta)\alpha_2(\eta\rho(r))$.

			\textsl{Step~5.} It remains to define the sequence $(T_i)_{i\in\N}$ of $(\sigalg_t)_{t\ge 0}$-optional times. Letting $T_{2k-1} := \check t_k$ and $T_{2k} := \hat t_k$, $k\in\N$, we see from Steps~2 through 4 that
			\begin{align*}
				& \Expec{V_{\sigma(t)}(x(t))\indic{\{t\in[T_{i-1}, T_i[\}\cap\{T_{i-1} < \infty\}}}\\
				& \qquad\qquad\qquad\qquad\le \beta(\norm{\xz}, t)\mx \gamma\bigl(\norm{d}_{\posR}\bigr),
			\end{align*}
			which proves the claim.\hfill{}$\square$

		\emph{Proof of Theorem~\ref{t:issmunif} (Sketch).}
		Fix $\nu\in\N$, and let $k' := \mu\left[\sum_{j\in\PSet}\frac{q_j}{\lambda_j}\left[1-\frac{1-\epower{-\lambda_j T}}{\lambda_j T}\right]\right]\big/\left[1-\sum_{j\in\PSet}\frac{\mu q_j\left(1-\epower{-\lambda_j T}\right)}{\lambda_j T}\right]$. In view of (Vd2), pointwise on $\bigl\{s\in[\tau_i, \tau_{i+1}[\bigr\}$, $i\in\N$,  and applying (Vd3) at $t = \tau_{i+1}$,
			\begin{align*}
				& V_{\sigma(\tau_{i+1})}(x(\tau_{i+1})) \le \mu V_{\sigma(\tau_i)}(x(\tau_i))\epower{-\lambda_{\sigma(\tau_i)}(\tau_{i+1}-\tau_i)}\\
				& \qquad\qquad + \frac{\mu \chi\bigl(\norm{d}_{\posR}\bigr)}{\lambda_{\sigma(\tau_i)}}\!\left(1-\epower{-\lambda_{\sigma(\tau_i)}(\tau_{i+1}-\tau_i)}\right).
			\end{align*}
			Iterating the above inequality from $i = 0$ through $i = \nu-1$, we get
			\begin{multline}
				V_{\sigma(\tau_\nu)}(x(\tau_\nu)) \le \mu^\nu V_{\sigma(0)}(\xz)\prod_{i = 0}^{\nu-1}\epower{-\lambda_{\sigma(\tau_i)}(\tau_{i+1}-\tau_i)}\\
				+ \mu^\nu\chi\bigl(\norm{d}_{\posR}\bigr) \sum_{i = 0}^{\nu-1} \frac{\mu^{-i}}{\lambda_{\sigma(\tau_i)}}\left(\prod_{j = i+1}^{\nu-1}\epower{-\lambda_{\sigma(\tau_j)}(\tau_{j+1}-\tau_j)}\right.\\
				\left.- \prod_{j = i}^{\nu-1}\epower{-\lambda_{\sigma(\tau_j)}(\tau_{j+1}-\tau_j)}\right).
			\label{e:Vtaunuunif}
			\end{multline}
			The expectation of the first term on the right-hand side of~\eqref{e:Vtaunuunif} can be evaluated as
			\begin{align}
				& \Expec{\mu^\nu V_{\sigma(0)}(\xz)\prod_{i = 0}^{\nu-1}\epower{-\lambda_{\sigma(\tau_i)}(\tau_{i+1}-\tau_i)}}\nonumber\\
				& \le \alpha_2(\norm{\xz})\left(\sum_{j\in\PSet}\frac{\mu q_j\left(1-\epower{-\lambda_j T}\right)}{\lambda_j T}\right)^\nu,
				\label{e:firstexpecunif}
			\end{align}
			by utilizing (Vd1) and (UH1)-(UH3). Also, from (UH3) we have
			\begin{align}
				& \Expec{\frac{\left(\prod_{j = i+1}^{\nu-1}\epower{-\lambda_{\sigma(\tau_j)}(\tau_{j+1}-\tau_j)} - \prod_{j = i}^{\nu-1}\epower{-\lambda_{\sigma(\tau_j)}(\tau_{j+1}-\tau_j)}\right)}{\lambda_{\sigma(\tau_i)}}}\nonumber\\
				& = \prod_{j=1+1}^{\nu-1}\Expec{\epower{-\lambda_{\sigma(\tau_{j+1})} S_{j+1}}}\Expec{\frac{1-\epower{-\lambda_{\sigma(\tau_i)} S_{i+1}}}{\lambda_{\sigma(\tau_i)}}}.
				\label{e:partial1unif}
			\end{align}
			Now for each $j\in\N$ we have
			\begin{align}
				\Expec{\epower{-\lambda_{\sigma(\tau_j)} S_{j+1}}} = \sum_{k\in\PSet} \frac{q_k\left(1-\epower{-\lambda_k T}\right)}{\lambda_k T},
				\label{e:partial2unif}
			\end{align}
			and for each $i\in\N$,
			\begin{align}
				\Expec{\frac{1-\epower{-\lambda_{\sigma(\tau_i)} S_{i+1}}}{\lambda_{\sigma(\tau_i)}}} = \sum_{k\in\PSet} \frac{q_k}{\lambda_k}\left(1-\frac{1-\epower{-\lambda_k T}}{\lambda_k T}\right).
				\label{e:partial3unif}
			\end{align}
			Substituting the right-hand sides of~\eqref{e:partial3unif} and~\eqref{e:partial2unif} back into~\eqref{e:partial1unif} and simplifying, we see that
			\begin{align}
				\Expec{V_{\sigma(\tau_\nu)}(x(\tau_\nu))} & \le \alpha_2(\norm{\xz})\left(\sum_{j\in\PSet}\frac{\mu q_j\left(1-\epower{-\lambda_j T}\right)}{\lambda_j T}\right)^\nu\nonumber\\
				& \quad + k'\chi\bigl(\norm{d}_{\posR}\bigr).
				\label{e:Vdestnuunif}
			\end{align}
			Now, letting $\gamma(r) := k'\chi(r)$ and $\beta(r, s) := \alpha_2(r)\eta^s$, where $\eta := \sum_{j\in\PSet}\frac{\mu q_j\left(1-\epower{-\lambda_j T}\right)}{\lambda_j T}$, an application of (Vd1) on the left-hand side of~\eqref{e:Vdestnuunif} immediately proves the assertion.
			\hfill$\square$

	\section{Control synthesis for \iss{} disturbance attenuation}
		\label{s:syn}
		We look at two different controller architectures, namely, one in which the controller is mode-dependent, and the other in which the controller is mode-independent. That is to say, in the first case, $u$ is a function of both the state $x$ and the switching signal $\sigma$, while in the second case $u$ is just a function of $x$.

		\subsection{Mode-dependent controllers} 
		Consider the affine-in-control switched system perturbed by a disturbance signal
		\begin{equation}
			\dot x = f_\sigma(x, d) + \sum_{i=1}^m g_{\sigma, i}(x) u_i,\qquad x(0) = \xz, \quad t\ge 0,
			\label{e:ssysdefcond}
		\end{equation}
		where $x\in\R^n$ is the state, $u_i$, $i=1, \ldots, m$, are the (scalar) control inputs, $f_j:\R^n\times\R^k\lra\R^n$ and $g_{j, i}:\R^n\lra\R^n$ are smooth maps for each $j\in\PSet$, $i\in\{1, \ldots, m\}$. Let $\CSet$ be the set where the control $u := [u_1, \ldots, u_m]\transp$ takes its values. For the moment we let $\CSet$ be a subset of $\R^m$ containing the origin. With a feedback control function $k_\sigma(x) := [u_{\sigma, 1}(x), \ldots, u_{\sigma, m}(x)]\transp$, the closed-loop system stands as
		\begin{equation}
			\dot x = f_\sigma(x, d) + \sum_{i=1}^m g_{\sigma, i}(x) k_{\sigma, i}(x),\qquad x(0) = \xz, \quad t\ge 0,
			\label{e:ssysdefcondcl}
		\end{equation}
		We let the switching signal $\sigma$ be a stochastic process as defined in~\secref{s:prelims}, and let $\xz\neq 0$.

\addtolength{\textheight}{0cm}

		Our goal is to choose a control function $k_\sigma$ so that~\eqref{e:ssysdefcondcl} satisfies some \iss{} in $\Lp 1$ estimate at switching instants. We shall appeal to our analysis results of~\secref{s:ana} and universal formulae for \iss{} disturbance attenuation to achieve this objective.

		Universal feedback control functions attaining \iss{} disturbance attenuation for nonlinear systems affected by disturbances and possessing control inputs were constructed in~\cite{ref:liberzonISSclf}. The results in that article rely on universal formulae for asymptotic feedback stabilization of nonlinear systems; applications include systems in which the control takes values in various restricted control sets, and a universal formula is available. In our illustrative result below we utilize off-the-shelf universal feedback control functions for \iss{} disturbance attenuation from~\cite{ref:liberzonISSclf}. The next proposition is a typical illustration of such a result.

			Let us define the map $\varphi:\R\times\R^m\lra\R$ given by
			\[
				\varphi(a, b) := 
				\begin{cases}
					\displaystyle{-\frac{a + \sqrt{a^2 + \norm{b}^4}}{\norm{b}^2}b} & \quad \text{if }b\neq 0,\\
					0 & \quad \text{otherwise},
				\end{cases}
			\]
			the function $\wt W_j(x) := \bigl[\LieD{g_{j, 1}}{V_j}(x), \ldots, \LieD{g_{j, m}}V_j(x)\bigr]$, and a map $\ol W_j:\R^n\lra\R$, with values chosen such that it is smooth away from $0$ and continuous at $0$, and
			\begin{align}
				&\max_{d\in\R^k}\left\{\frac{\partial V_j}{\partial x}(x)f_j(x, d) - \chi(\norm d)\right\} + \lambda_j V_j(x) \le \ol W_j(x)\nonumber\\
				&\quad \le \max_{d\in\R^k}\left\{\frac{\partial V_j}{\partial x}(x)f_j(x, d) - \chi(\norm d)\right\} + 2\lambda_j V_j(x)
				\label{e:longdef}
			\end{align}
			for all $x\in\R^n$, $j\in\PSet$.
		\begin{proposition}
			Consider the system~\eqref{e:ssysdefcond} with $\CSet = \R^m$. Suppose that $\sigma$ belongs to class UH, and
			\begin{enumerate}[{\rm (Cd1)}]
				\item (Vd1) of Assumption~\ref{a:Vi} holds;
				\item (Vd3) of Assumption~\ref{a:Vi} holds;
				\item $\therex \alpha, \chi\in\ClassKinfty$, $\therex \lambda_j\in\Lambda = \R$, $j\in\PSet$, such that $\fa x\in\R^n\setmin\{0\}$, $\fa d\in\R^k$ and $\fa j\in\PSet$ we have
					\begin{align*}
						& \inf_{u\in\CSet}\left\{\frac{\partial V_j}{\partial x}(x)f_j(x, d) + 3\lambda_j V_j(x)\right.\\
						& \qquad\quad\quad\quad \left.+ \sum_{i=1}^m \LieD{g_{j, i}}V_j(x) u_i\right\} \le \chi(\norm d);
					\end{align*}
				\item $\fa \eps > 0 \;\therex \delta > 0$ such that if $x(\neq 0)$ satisfies $\norm{x} < \delta$, then $\therex u\in\R^m,\; \norm u < \eps$, such that $\fa j\in\PSet$
					\begin{align*}
						& \max_{d\in\R^k}\left\{\frac{\partial V_j}{\partial x}(x)f_j(x, d) - \chi(\norm d)\right\}\\
						& \qquad\quad + \sum_{i = 1}^m \LieD{g_{j, i}}{V_j}(x) u_i \le -\lambda_j V_j(x);
					\end{align*}
				\item (U3) of Theorem~\ref{t:issmunif} hold.
			\end{enumerate}
			Then under the feedback control function
			\begin{equation}
				k_\sigma(x) = \varphi\left(\ol W_\sigma(x), \wt W_\sigma\transp(x)\right)
				\label{e:univformulaiss}
			\end{equation}
			the system~\eqref{e:ssysdefcondcl} satisfies an \iss{} in $\Lp 1$ estimate at switching instants, 
			\label{p:issmclsigmadep}
		\end{proposition}

		The proof relies heavily on the proof of~\cite[Theorem~3]{ref:liberzonISSclf}, see~\cite[Chapter~3]{ref:myphdthesis} for details.

		\subsection{Mode-independent controllers.}
		Consider the affine in control switched system~\eqref{e:ssysdefcond}. Let $k(x) = [k_1(x), \ldots, k_m(x)]\transp$ be a feedback control function, with which the closed-loop system stands as
		\begin{equation}
			\dot x = f_\sigma(x, d) + \sum_{i=1}^m g_{\sigma, i}(x)k_i(x),\;\;\; x(0) = \xz,\;\; t\ge 0.
			\label{e:ssysdefcondclunob}
		\end{equation}
		We let the switching signal $\sigma$ be a stochastic process as defined in~\secref{s:prelims}, and let $\xz\neq 0$.

		Our objective is to choose a control function $k$ such that~\eqref{e:ssysdefcondclunob} satisfies an \iss{} in $\Lp 1$ estimate at switching instants, for some class-$\ClassKinfty$ function $\alpha$.

		\begin{proposition}
			Consider the system~\eqref{e:ssysdefcond} with $\CSet = \R^m$. Suppose that $\sigma$ belongs to class UH, and
			\begin{enumerate}[{\rm (CUd1)}]
				\item (Vd1) and (Vd3) of Assumption~\ref{a:Vi} holds;
				\item there exists a control function $k:\R^n\lra\CSet$, such that $\displaystyle{\frac{\partial V_i}{\partial x}(x)\bigl(f_i(x, d) + g_i(x)k(x)\bigr) \le -\lambda_i V_i(x)} + \chi(\norm d)$ for every $i\in\PSet$, $x\in\R^n$;
				\item $\displaystyle{\sum_{i\in\PSet}\frac{\mu q_i\bigl(1-\epower{-\lambda_i T}\bigr)}{\lambda_i T} < 1}$.
			\end{enumerate}
			Then under $k$ the system~\eqref{e:ssysdefcond} satisfies an \iss{} in $\Lp 1$ estimate at switching instants.
			\label{p:sigmaunifissonobcon}
		\end{proposition}

		The assertion follows immediately by first observing that the closed-loop system is~\eqref{e:ssysdefcondclunob}, and then applying Theorem~\ref{t:issmunif} to~\eqref{e:ssysdefcondclunob}.

	\section{Remarks on Markovian Switching Signals}
		\label{s:concl}
		We have established sufficient conditions for different \iss{}-type properties of the randomly switched system~\eqref{e:ssysdef} under different classes of switching signals. Let us reiterate that for switching signals of class G and UH considered here, it is difficult to write an infinitesimal generator, since there is either too little information about the parameters of the switching signal, or a strong dependence on its past history. For switching signals coming from continuous-time Markov chains it is possible to employ the infinitesimal (or extended) generator to derive conditions for stability. We get stronger bounds by this route, as shown in~\cite[Chapter~3]{ref:myphdthesis}. Indeed, we have
		\begin{theorem}
			Consider the system~\eqref{e:ssysdef}, and suppose that $\sigma$ is a continuous-time Markov chain with generator matrix $Q = [q_{i, j}]_{\cardP\times\cardP}$. Moreover, suppose that there exist functions $V:\PSet\times\R^n\lra\posR$, $V(i, \cdot)$ is continuously differentiable for each $i$, $\alpha_1,\alpha_2, \rho\in\ClassKinfty$, and a constant $\lambda_\circ > 0$, such that
			\begin{itemize}
				\item $\alpha_1(\norm x) \le V(i, x) \le \alpha_2(\norm x)$,
				\item $\mc LV(i, x) \le -\lambda_\circ V(i, x)$ whenever $\norm{x} \ge \rho(\norm d)$.
			\end{itemize}
			Then the inequality in~\eqref{e:ideal} holds for some $\beta\in\ClassKL$ and some $\alpha, \gamma\in\ClassKinfty$.
		\end{theorem}

		For definitions of Markov chains, (local) martingales, and martingale problems, see, e.g.,~\cite{ref:raoStochProcGenTheo}. The operator $\mc L$ is defined in terms of an appropriate martingale problem as follows. Let $h:\PSet\times\R^n\lra\R$ be a function such that there exists a measurable function $\wt h:\PSet\times\R^n\lra\R$ such that the process
		\[
		\left(h(\sigma(t), x(t)) - h(\sigma_0, \xz) - \int_0^t \wt h(\sigma(s), x(s))\drv s\right)_{t\ge 0}
		\]
		is a mean-zero $(\sigalg_t)_{t\ge 0}$-local martingale. We define $\mc Lh(i, x) := \wt h(i, x)$, where $\mc L$ is the extended generator~\cite{ref:davisMarkovModelsOptbk} corresponding to the Markov process $(\sigma(t), x(t))_{t\ge 0}$. Of course, finding the class of functions $h$ for which such a $\wt h$ exists is a nontrivial matter, but it is usually not difficult to find a subclass. Often the operator $\mc L$ is defined in terms of a differentiation operation, namely,
		\[
		\mc Lh(i, x) = \lim_{h\downarrow 0}\frac{\Expec{h(\sigma(t+h), x(t+h))\big|A_t} - h(i, x)\bigr)}{h},
		\]
		where $A_t := \bigl\{(\sigma(t), x(t)) = (i, x)\bigr\}$, and $h:\PSet\times\R^n\lra\R$ is a function that is pointwise continuously differentiable on the set $\PSet$.

		A similar approach relying on the solution to appropriate martingale problems can be adopted if $\sigma$ is a general marked point process~\cite{ref:bremaudPointProc} with suitable stochastic jump intensities, and will be reported elsewhere. Another interesting direction of work concerns establishing \iss{}-type estimates ``in probability'' of~\eqref{e:ssysdef}, such as those formulated in~\cite{ref:krsticdengbk, ref:smallgainwrong}. 

\bibliographystyle{IEEEtran}
\bibliography{refspec}

\begin{thebibliography}{10}
\providecommand{\url}[1]{#1}
\csname url@rmstyle\endcsname
\providecommand{\newblock}{\relax}
\providecommand{\bibinfo}[2]{#2}
\providecommand\BIBentrySTDinterwordspacing{\spaceskip=0pt\relax}
\providecommand\BIBentryALTinterwordstretchfactor{4}
\providecommand\BIBentryALTinterwordspacing{\spaceskip=\fontdimen2\font plus
\BIBentryALTinterwordstretchfactor\fontdimen3\font minus
  \fontdimen4\font\relax}
\providecommand\BIBforeignlanguage[2]{{%
\expandafter\ifx\csname l@#1\endcsname\relax
\typeout{** WARNING: IEEEtran.bst: No hyphenation pattern has been}%
\typeout{** loaded for the language `#1'. Using the pattern for}%
\typeout{** the default language instead.}%
\else
\language=\csname l@#1\endcsname
\fi
#2}}

\bibitem{ref:sontagISS}
E.~D. Sontag, ``Smooth stabilization implies coprime factorization,''
  \emph{IEEE Transactions on Automatic Control}, vol.~34, pp. 435--443, 1989.

\bibitem{ref:sontagISSnote}
------, ``Input-to-state stability: basic concepts and results,'' in
  \emph{Nonlinear and Optimal Control Theory}, P.~Nistri and G.~Stefani,
  Eds.\hskip 1em plus 0.5em minus 0.4em\relax Springer-Verlag, 2006, pp.
  163--220.

\bibitem{ref:tsiniasExpISS}
J.~Tsinias, ``The concept of `exponential input to state stability' for
  stochastic systems and applications to feedback stabilization,''
  \emph{Systems {\&} Control Letters}, vol.~36, pp. 221--229, 1999.

\bibitem{ref:liberzonbk}
D.~Liberzon, \emph{Switching in Systems and Control}, ser. Systems \& Control:
  Foundations \& Applications.\hskip 1em plus 0.5em minus 0.4em\relax Boston:
  Birkh{\"a}user, 2003.

\bibitem{ref:linhAutomatica06}
L.~Vu, D.~Chatterjee, and D.~Liberzon, ``Input-to-state stability of switched
  systems and switching adaptive control,'' \emph{Automatica}, vol.~43, no.~4,
  pp. 639--646, April 2007.

\bibitem{ref:davisMarkovModelsOptbk}
M.~H.~A. Davis, \emph{Markov Models and Optimization}.\hskip 1em plus 0.5em
  minus 0.4em\relax London: Chapman \& Hall, 1993.

\bibitem{ref:ransw}
\BIBentryALTinterwordspacing
D.~Chatterjee and D.~Liberzon, ``Stability analysis of randomly switched
  systems,'' To appear in {IEEE} {T}ransactions on {A}utomatic {C}ontrol, 2007.
  [Online]. Available:
  \url{http://decision.csl.uiuc.edu/~liberzon/publications.html}
\BIBentrySTDinterwordspacing

\bibitem{ref:stabranswcorr}
\BIBentryALTinterwordspacing
------, ``Stabilizing randomly switched systems,'' submitted for publication,
  2006. [Online]. Available:
  \url{http://decision.csl.uiuc.edu/~liberzon/publications.html}
\BIBentrySTDinterwordspacing

\bibitem{ref:liberzonISSclf}
D.~Liberzon, E.~D. Sontag, and Y.~Wang, ``Universal construction of feedback
  laws achieving {ISS} and integral-{ISS} disturbance attenuation,''
  \emph{Systems \& Control Letters}, vol.~46, pp. 111--127, 2002.

\bibitem{ref:raoStochProcGenTheo}
M.~M. Rao, \emph{Stochastic {P}rocesses: {G}eneral {T}heory}, ser. Mathematics
  and Its Applications.\hskip 1em plus 0.5em minus 0.4em\relax Dordrecht:
  Kluwer Academic Publishers, 1995, vol. 342.

\bibitem{ref:bremaudPointProc}
P.~Br{\'e}maud, \emph{Point Processes and Queues, Martingale Dynamics}, ser.
  Springer Series in Statistics.\hskip 1em plus 0.5em minus 0.4em\relax New
  York - Berlin: Springer-Verlag, 1981.

\bibitem{ref:borovkovbk}
A.~A. Borovkov, \emph{Probability Theory}.\hskip 1em plus 0.5em minus
  0.4em\relax Amsterdam: Gordon \& Breach Publishing Group, 1999.

\bibitem{ref:myphdthesis}
\BIBentryALTinterwordspacing
D.~Chatterjee, ``Studies on stability and stabilization of randomly switched
  systems,'' Ph.D. dissertation, University of Illinois at Urbana-Champaign,
  2007. [Online]. Available:
  \url{http://decision.csl.uiuc.edu/~liberzon/collaborators.html}
\BIBentrySTDinterwordspacing

\bibitem{ref:pralywangISS}
L.~Praly and Y.~Wang, ``Stabilization in spite of matched unmodelled dynamics
  and an equivalent definition of input-to-state stability,'' \emph{Mathematics
  of Control, Signals and Systems}, vol.~9, pp. 1--33, 1996.

\bibitem{ref:raoProbTheo}
M.~M. Rao and R.~J. Swift, \emph{Probability {T}heory with {A}pplications},
  2nd~ed., ser. Mathematics and Its Applications.\hskip 1em plus 0.5em minus
  0.4em\relax Springer-Verlag, 2006, vol. 582.

\bibitem{ref:krsticdengbk}
M.~Krsti{\'c} and H.~Deng, \emph{Stabilization of Nonlinear Uncertain
  Systems}.\hskip 1em plus 0.5em minus 0.4em\relax Springer-Verlag, 1998.

\bibitem{ref:smallgainwrong}
Z.-J. Wu, X.-J. Xie, and S.-Y. Zhang, ``Adaptive backstepping controller design
  using stochastic small-gain theorem,'' \emph{Automatica}, vol.~43, no.~4, pp.
  608--620, 2007.

\end{thebibliography}

\end{document}